\documentclass[twocolumn, 10pt,pre]{revtex4-1}

\usepackage{algcompatible}
\usepackage{newfloat}
\usepackage{amsmath}
\usepackage{amssymb}
\usepackage{amsthm}
\usepackage{graphicx}
\usepackage{subfig}
\usepackage{latexsym}
\usepackage{epsfig}
\usepackage{color}

\DeclareFloatingEnvironment[
    fileext=loa,
    listname=List of Algorithms,
    name=ALGORITHM,
    placement=tbhp,
]{algorithm}

\def\ind{\bf{1}}

\def\u2{U_{2}}
\def\u1{U_{1}}

\newtheorem{teo}{Theorem}[section]

\begin{document}

\title{On minimum correlation in construction of multivariate distributions}
\author{Vanja M. Dukic}
\affiliation{Department of Applied Mathematics, University of Colorado, Boulder, Colorado 80309, USA}
\author{Nevena Mari\'c}
\affiliation{Department of Mathematics, University of Missouri-St. Louis, St. Louis, Missouri 63121, USA}

\begin{abstract}
In this paper we present an algorithm for exact generation of multivariate samples with pre-specified marginal distributions and  a given correlation matrix, based on a mixture of Fr\'echet-Hoeffding bounds and marginal products.  The algorithm can accommodate any among the theoretically possible correlation coefficients, and explicitly provides a connection between simulation and the minimum correlation attainable for different distribution families. We calculate the minimum correlations  in several common distributional examples, including in some that have not been looked at before. As an illustration, we provide the details and results of implementing the algorithm for generating three-dimensional negatively and positively correlated Beta random variables, making it the only non-copula algorithm for correlated Beta simulation in dimensions greater than two.  This work has potential for impact in a variety of fields where simulation of multivariate stochastic components is desired.
\end{abstract}

\maketitle

\section{Introduction}

The original question of generating multivariate probability distributions with pre-specified margins has a long history, dating in part back to the work of E. Wigner \cite{Wigner1932} on thermodynamic equilibria.  The general form of the problem was studied by M. Fr\'echet \cite{frechet} and V. Hoeffding \cite{hoeffding}   in a body of work which grew out of the problem originally  posed by P. L\'evy \cite{levy}. Today, this work falls under the scope of Fr\'echet-Hoeffding classes. An excellent  overview of the developments in this field can be found in  Dall'Aglio {\it et al.} \cite{dall},  R\"uschendorf  {\it et al.} \cite{Ruschendorf1996}, and Conway \cite{Conway1979}. 

Today, algorithms for generation of correlated random variables with pre-specified marginal distributions play an important role in simulation of stochastic processes, and hybrid deterministic-stochastic systems. Such algorithms are encountered in a variety of fields, for example: statistics and applied probability \cite{Craiu2001,Craiu2005,Schmei1982,Rosenfeld2008}, finance \cite{Lawrance1981}, environmental science \cite{Izawa1965}, physics \cite{Smith1981}, engineering \cite{Lampard1968}, and ecology where "demographic" or "weather" stochasticity is an increasingly more relevant component of species dynamics  \cite{Tadeu2008}. Much of the development of these algorithms  has so far relied on coupling ideas -- or antithetic coupling for negatively-correlated variables \cite{Thorisson1991} -- and copula-based methods \cite{Series2010,Nelsen2006a}.  

{ Copula methods, in particular, have recently become widely used in generation of samples from 
multivariate distributions with pre-specified marginals and a dependency function \cite{Nelsen2006}. Copula methodology relies on the results of Sklar \cite{Series2010} 
who proved that a multivariate distribution can be characterized  (uniquely in the case of continuous distributions)  
 by the set of its marginal distributions and a "copula" function which describes the dependence between the 
 components.
 The dependence among the original variables  is then translated, via  the copula function, into the dependence on the 
 scale of uniform random variables. Consequently, the entire desired multivariate distribution is obtained via a 
transformation of these correlated uniform variables.  Unfortunately, the correlation is not preserved 
under these transformations, and the sampling is not exact.}

In this paper we  present a novel alternative algorithm that  generates   {\it exact \,} multivariate samples
 with pre-specified marginals and a given  correlation matrix. We note that
 specifying marginal distributions and a correlation matrix is in general not enough to completely determine 
the entire multivariate distribution. Nonetheless, specifying marginal distributions and a set of linear
 relationships (through a set of correlation coefficients) among the random variables is frequently done, 
perhaps due to the strong intuitiveness of the linear relationship.  

The algorithm can be used  for generating a realization of a set of random variables 
$X_1,X_2,... X_t$ such that each variable $X_i$ has a specified marginal distribution $F_i$, 
and such that each pairwise correlation coefficient $cor(X_i,X_j)$ equals some set value $\rho_{ij}$. 
We take the correlation coefficient between two random variables $X_i$ and $X_j$ to  be defined as
$$cor(X_i,X_j) = \frac{cov(X_i,X_j)}{\sigma_{X_i} \sigma_{X_j}},$$
where $cov(X_i,X_j)$ stands for covariance between $X_i$ and $X_j$.

The paper is organized as follows: Section~\ref{algorithm} introduces the basic idea of the algorithm, 
which allows for fast simulation and can accommodate any among the theoretically plausible  correlation ranges. 
We discuss its implementation and performance, and present detailed examples for several bivariate distribution 
families (Uniform, Arcsine, Weibull, Exponential, Erlang, Beta and Gaussian) in Section~\ref{examples}. We also 
calculate the minimal correlations (or {\it maximum negative} correlations as they are called in Kotz {\it et al.} 
\cite{Kotz2000}) for these distributions, including some that have never been obtained before.  In Section~\ref{multi} 
we present a multivariate  extension of the algorithm.  Finally, Section~\ref{Conclusions} concludes the paper with a 
brief summary of benefits and limitations of the proposed approach.

\section{The algorithm for generating bivariate samples from prescribed marginals with specified correlation}
\label{algorithm}

Perhaps not surprisingly, over the years the problem of generating bivariate distributions with fixed marginals and a specified correlation coefficient has gotten more attention from the simulation  communities (e.g \cite{Hill1994}), than from the probability community \cite{Devroye1996}. A notable exception is a work by Johnson and Tenenbein \cite{Johnson1981}, who provide a bivariate generation algorithm based on a different method than the one presented in this paper. Moreover, bivariate distributions have  been studied mostly for marginals with common distributional form, such as the Normal, Exponential, or Gamma (see e.g. \cite{Kotz2000}).  The case of exponential marginals is particularly well studied; however,  providing a constructive algorithm that can produce all theoretically possible correlation coefficients still proved to be a significant theoretical advance \cite{Bladt2010}.  For general marginals, Michael and Schucany  \cite{Michael2002} introduced a  mixture-density based (hierarchical) approach, although this algorithm relies on finding a feasible mixture density for each example, which need not always be straightforward.  

Trivariate reduction methods were introduced initially for construction of  dependent Gamma random variable pairs (see for example 
\cite{Arnold1967,Schmei1982,Lai1995}), as an alternative to once computationally costly {\it distribution inversion} based methods. 
The trivariate reduction idea relies on the use of three independent variables in order to obtain one pair of  correlated variables.
However, these methods are limited to additive families of distributions, like Gamma or Poisson. 
The algorithm we present here  will not have that limitation. 

The algorithm in this paper is a hybrid version of the trivariate reduction method, as it  relies on three uniformly distributed random variables to produce a pair, but it is not inversion-free.  It is based on the following reasoning: with a certain probability we use the same source of randomness in the construction of the pair, and two independent sources otherwise. The probability used to determine which source is  used will be closely related to the correlation coefficient. 

To set notation, let $F$ and $G$ be cumulative distribution functions (cdfs) with finite positive variances, and let $X$ and $Y$ be random
variables with distributions $F$ and $G$ respectively, $X\sim F$ and $Y \sim G$. The first question to be asked is whether 
 {\em any correlation} $\rho \in [-1,1]$ can be attained for the pair $(X,Y)$. The answer to that question is negative, 
and dates back to the work of Hoeffding \cite{hoeffding} and Fr\'echet \cite{frechet}, where  the concept of {\em extremal 
distributions} was originally introduced: if we let $\Pi(F,G)$ be the set of all bivariate cdfs having $F$ and $G$ as marginals,
 then among the elements of $\Pi(F,G)$, there are cdfs $H^\ast$ and $H_\ast$  which have maximum and minimum correlation coefficient 
($\rho^\ast$ and $\rho_\ast$), respectively.  Such extremal distributions are also called (upper and lower) {\em Fr\'echet-Hoeffding bounds}.  
They were later characterized by Whitt \cite{Whitt1976} who provides the following two equivalent statements.
\begin{teo}[Hoeffding (1940)]\label{whitt1}
For any F and G with finite positive variances,
\begin{eqnarray*}
&&  H^\ast (x,y)=\min\{F(x), G(y)\} \\
&& H_\ast(x,y)= \max\{0,[F(x)+G(y)-1]\}, 
\end{eqnarray*}
for all $(x,y) \in \mathbb{R}^2.$
\end{teo}
\begin{teo}[Hoeffding (1940); Whitt (1976)] \label{whitt2}
For any F and G with finite positive variances
\begin{eqnarray*}
&& (F^{-1}(U), G^{-1}(U)) ~\mbox{has cdf}~ H^\ast \\
&&\mbox{and}~ (F^{-1}(U), G^{-1}(1-U))~ \mbox{has  cdf}~ H_\ast 
\end{eqnarray*}
where $U$ has uniform distribution on $[0,1]$ and $F^{-1}$ and $G^{-1}$ are inverse distribution functions, defined as
 $F^{-1}(y)=\inf\{x: F(x) \geq y \}$, and $G^{-1}(y)=\inf\{x: G(x) \geq y \}$, respectively.
\end{teo}
Fr\'echet \cite{frechet} suggested that any system of bivariate distributions with specified marginals $F$ and $G$ should include $H_\ast$ and $H^\ast$ as limiting cases \cite{Kotz2000}.  The crux of our algorithm is precisely in this reasoning, as we construct multivariate distributions through careful blending of Fr\'echet-Hoeffding bounds and marginal products. This blending, although clearly apparent in the bivariate case, becomes less obvious in dimensions greater than two and in the presence of negative correlations.

\subsection{The Basics: Bivariate Algorithm}
Suppose $F$ and $G$ are desired marginal distributions, with finite positive variances. 
Then the main bivariate problem can be stated as follows:
Construct $X$ and $Y$ such that $X \sim F$, $Y \sim G$, and correlation $cor(X,Y)=\rho$. Here, 
$\rho \in [\rho_\ast,\rho^\ast]$, where $\rho_\ast$ and $\rho^\ast$ are minimum and maximum theoretically 
possible correlation coefficients, respectively.

Let $m_F, \sigma_F$ and $m_G, \sigma_G$ be the first moments and standard deviations corresponding to $F$ and $G$,  respectively. 
Let $\phi$ be an algorithm such that $\phi(U) \sim F$, and  let $\psi$ be an algorithm such that $\psi(U) \sim G$, where $U$ is a uniformly distributed random variable on $[0,1]$. (It can be assumed, although it is not necessary, that $\phi=F^{-1}$ and $\psi=G^{-1}$). Let $V$ also be a uniform random variable  on $[0,1]$ and define
\begin{eqnarray}\label{defc}
c^{\phi,\psi}(U,V)= \frac { E[\phi(U)\psi(V)] -m_F m_G  }{\sigma_F \sigma_G},
\end{eqnarray}
where $E[\cdot]$  is used to denote the expected value of a random variable.
To simplify notation, we will denote $c^{\phi,\phi}$ as $c^{\phi}$. 
Also, observe that $c^{\phi}(U,U) =(E(X^2) - m_{F}^2)/\sigma^2_{F} = 1$, for $X = \phi(U) \sim F$.  

The following construction, presented as Algorithm 1 below, will yield a pair of variables, $(X,Y)$, such that  $X \sim F$, $Y \sim G$ and $cor(X,Y)=\rho$ for  $\rho \in [\rho_\ast,\rho^\ast]$.

\begin{algorithm} [h]
\hrule
\caption{Construction of two random variables with prescribed marginal distributions $F, G$ 
and correlation coefficient $\rho$\\[-1ex]}  \hrule \label{alg2}
\begin{algorithmic}[1]
\STATE {\bf sample} $U,V,W \sim U(0,1)$, independently
\STATE {\bf let} $X=\phi(U)$
\IF {$\rho> 0$} \label{volta}
\STATE  {\bf let} $U'=  U$.
\ELSE
\STATE {\bf let} $U'= 1- U$.
\ENDIF
\IF {$W < \rho/ c^{\phi,\psi}(U,U')$}
\STATE {\bf let} $Y=\psi(U')$ \label{state1}
\ELSE
\STATE {\bf let} $Y=\psi(V)$ \label{state2}
\ENDIF
\STATE RETURN $(X,Y)$ 
\vspace{.1in}
\hrule
\end{algorithmic}
\end{algorithm}

\begin{teo}\label{propo1}
If  $(X,Y)$ is generated by Algorithm~\ref{alg2}, then 
\begin{enumerate}
\item[(a)] $X \sim F$, $Y \sim G$ and $cor(X,Y)=\rho$, if $\rho/c^{\phi,\psi}(U,U')\leq 1$.
\item[(b)] If $\phi=F^{-1}, \psi=G^{-1}$ then
\begin{enumerate} 
\item[i)] $(X,Y)$ has the joint distribution $H(x,y)$, where:
\begin{eqnarray*}
\mbox{for \,}&\rho \geq 0:&H(x,y)= \frac{\rho}{\rho^\ast} H^\ast(x,y) + (1-\frac{\rho}{\rho^\ast}) F(x) G(y);\\ 
\mbox{for \,}&\rho \leq 0:&H(x,y)= \frac{\rho}{\rho_\ast} H_\ast(x,y) + (1-\frac{\rho}{\rho_\ast}) F(x) G(y).
\end{eqnarray*}
\item[ii)]  Algorithm~\ref{alg2} is applicable for all $\rho_\ast \leq \rho \leq \rho^\ast$.
\end{enumerate}
\end{enumerate}
\end{teo}

\noindent {\bf Proof}: \\
{\bf (a)} By construction, $X \sim F$, $Y\sim G$. Using ${\ind}(\cdot)$ to denote the indicator function, $E[\cdot]$ to denote expected value of a random quantity, and $c$ in the place of $\rho/c^{\phi,\psi}(U,U')$, we have:
\begin{eqnarray*}
\lefteqn{E[XY] = E[XY{\ind}(W<c)] + E[XY{\ind}(W> c)]}\\
&&= E[\phi(U)\psi(U'){\ind}(W < c)] + E[\phi(U)\psi(V){\ind}(W> c)] \nonumber\\
&&\mbox{(U,V, W independent) }\\
&&=P(W < c)~ E[\phi(U)\psi(U')] + P(W > c)~ E\phi(U)~ E\psi(V) \nonumber\\
&&=c E[\phi(U)\psi(U')] +(1-c)m_F m_G  \nonumber \\
&&=c(E[\phi(U)\psi(U')] -m_F m_G) +m_F m_G.
\end{eqnarray*}
Then, from  (\ref{defc}),  it follows that $ cor(X,Y) = \rho.$
Observe also that when $\phi$ and $\psi$ are non-decreasing functions,  $cov(\phi(U),\psi(U))$ is always positive, while  $cov(\phi(U),\psi(1-U))$  is always negative. This can be verified easily using a coupling argument  as in  \cite{Thorisson1991}. Inverse distribution functions $\phi=F^{-1}, \psi=G^{-1}$ are of course non-decreasing.\\

{\bf (b)} For positive $\rho$,   Algorithm 1 produces $(F^{-1}(U),G^{-1}(U))$ with probability $\rho/c^{F^{-1},G^{-1}}(U,U)$. By Theorem \ref{whitt2}, the pair $(F^{-1}(U),G^{-1}(U))$   has the cdf $H^\ast (x,y)$. With probability $1-\rho/c^{F^{-1},G^{-1}}(U,U)$, the outcome of Algorithm 1 is a pair of two independent variables $(F^{-1}(U),G^{-1}(V))$, with the cdf that is a product of the marginal cdfs $F(x)G(y)$. The argument works analogously for negative values of $\rho$. \\

When $\rho>0$, 
by Theorem \ref{whitt2}, the maximum correlation between $F$ and $G$ is attained with the coupling $(F^{-1}(U), G^{-1}(U))$, so  that
\begin{equation}
\rho^\ast = \frac{ E[F^{-1}(U)G^{-1}(U)] -m_F m_G  }{\sigma_F \sigma_G}= c^{F^{-1},G^{-1}}(U,U) . 
\end{equation}
It follows that $\rho/c^{F^{-1},G^{-1}}(U,U) \leq 1 \Leftrightarrow \rho \leq \rho^\ast$.\\

When $\rho <0$, we again have, by  Theorem \ref{whitt2}, that the minimal correlation between $F$ and $G$ is attained with the coupling $(F^{-1}(U), G^{-1}(1-U))$, and that then
\begin{equation*}
\rho_\ast = \frac{ E[F^{-1}(U)G^{-1}(1-U)] -m_F m_G  }{\sigma_F \sigma_G} =
\end{equation*}\\[-6ex]
\begin{equation}
c^{F^{-1},G^{-1}}(U,1-U). \label{mincor}
\end{equation}
In this case it follows that $\rho/ c^{F^{-1},G^{-1}}(U,1-U) \leq 1 \Leftrightarrow \rho \geq  \rho_\ast.$
As noted in the last paragraph of the proof of part $(a)$, since
 $\phi=F^{-1}, \psi=G^{-1}$, we have that $\rho^\ast >0$ and $\rho_\ast <0$. Therefore the algorithm works for the entire range of possible correlations between $F$ and $G$. \hfill $\blacksquare$\\

\noindent{\bf{Remark}}:\\
\noindent Note that we allow a possibility that $\phi$ and $\psi$ are not inverse distribution functions, because the main idea of the algorithm is applicable  to transformations that are not inverse distribution functions.  We will present Algorithm 4 in  Section 3 as an example of using such transformations.

\subsection{Examples: Finding Minimum Correlations}\label{examples}

We now illustrate the implementation of Algorithm 1 using several common distributions as examples. We will
assume identical marginal distributions, $F=G$, and that  $\phi$  denotes an inverse distribution function in each case below.
Since the correlation coefficient is not  preserved under inverse distribution transformation -- namely, in general $cor(F^{-1}(U),F^{-1}(V)) \neq cor(U,V)$ --  the range of possible correlations for any individual distribution  has to be derived separately. Once the range of feasible correlations is known, application of Algorithm 1 is very simple and requires only few lines of code. 

It should be noted that determination of minimum (and maximum) possible correlation among two distributions has had a theoretical value in its own right. At the same time it is also of practical value, as knowing the maximum and minimum correlations  allows us to place the correlation estimates in perspective, which is  of great importance in empirical data analysis.   Moran \cite{Moran1967a} showed that only symmetric bivariate distributions for which there exist $\eta_0$ and $\eta_1$ such that $\eta_0 + \eta_1 Y$ has the same distribution as $X$  allow $\rho$ to take any value in the entire range $[-1,1]$. Some ranges for correlation coefficients for  bivariate distributional families are provided in \cite{Conway1979}; however, many ranges  still remain to be computed.
 
When the marginals are equal, maximum correlation $\rho^\ast =1$ since $cor(X,X)=1$, and only $\rho_\ast$ has to be determined.
We present briefly  several examples and derive the minimum correlation for each case. The first two cases, 
the Uniform and Arcsine, easily follow from Moran \cite{Moran1967a}, so we show them only as illustrations. 

\begin{itemize}
\item {\bf Uniform.} In the case of Uniform distribution on [0,1], we have that $\phi(U)=U$ and
\begin{eqnarray*}
\rho_\ast=\frac{E(U(1-U))-[E(U)]^2}{Var(U)}=\frac{1/6-1/4}{1/12} = -1.
\end{eqnarray*}

\item  {\bf Arcsine.} In the case of the Arcsine distribution with density $1/(\pi \sqrt{1-x^2})$ on [-1,1], $\phi(U)=\cos(\pi U)$ \cite{Devroye1996}. 
As the mean of this density is $0$ and variance $1/2$, it follows that
\begin{eqnarray*}
&\int_0^{1} \cos(\pi x) \cos\pi(1-x)dx= \int_0^1 \cos^2(\pi x) dx=&\\
& \frac{1}{2}\int_0^1 1+\cos(2\pi x) dx = \frac{1}{2}.&
\end{eqnarray*}
 From here $\rho_\ast=-1$, and that the algorithm is applicable for all $\rho \in [-1,1]$.

\item {\bf Exponential.} 
 If the variables are exponentially distributed with mean 1 { (density $e^{-x}$)} then  $\phi(U)=-log(U)$ and $\rho_\ast = 1- \pi^2/6  \approx -0.6449$, since
\begin{eqnarray*}
&E[\phi(U)\phi(1-U)] = \int_0^1 \log(x) \log(1-x) dx =& \\
&2 -\frac{1}{6} \pi^2.&
\end{eqnarray*}
The above integral can be solved using Maclaurin series representation of $\log(x)$, using either double or single series, and we present this proof in the Appendix.

Consequently, since $S \sim Exp(\lambda)$ (for $\lambda >0$) can be obtained as $\lambda  T$ where $T\sim Exp(1)$, the same range of attainable correlation $\rho \geq \rho_\ast= 1- \pi^2/6$  is valid for any choice of marginal exponential distribution. 

It is worth noting that many different bivariate exponential distribution algorithms have been studied, including a classic example by Gumbel \cite{Gumbel1960} and many others mentioned in \cite{Kotz2000}. Another recent construction of bivariate exponential distribution that allows an arbitrary positive or negative correlation coefficient has been introduced  by \cite{Bladt2010}. They use an elegant concept of  multivariate phase-type distributions, and provide a constructive algorithm that achieves minimum correlation $\rho_\ast$ through a limit of a sequence. 

\item {\bf Erlang.} As a $Gamma(n,\lambda)$ distribution where $n$ is an integer, an Erlang$(n,\lambda)$ random variable 
{(density $ x^{n-1} e^{-x/{\lambda}} / ((n-1)!\lambda^n) $)}
can be obtained as a sum of $n$ independent Exponential random variables  with parameter $\lambda$. 
Let $(S_1^1, S_1^2),(S_2^1, S_2^2),..., (S_n^1, S_n^2)$ be $n$ independent outputs of  Algorithm~1, where each variable in the pair has an exponential marginal distribution (with parameter $\lambda$), and where  $cor(S_i^1, S_i^2)=\rho$, for $i=1,...,n$. (Notice also that for $j \neq i$, $S_i^1$ and $S_j^1$ are independent, as are  $S_i^1$, and $S_j^2$.)
Let $X=S_1^1+...+S_n^1$ and $Y=S_1^2+...+S_n^2$. Then $X,Y \sim Gamma(n,\lambda)$ and $cor(X,Y)=\rho.$ It follows that the minimal possible correlation of $X$ and $Y$ is $1- \pi^2/6  \approx -0.6449$.

\item  {\bf Weibull. } The Weibull distribution with  density $k x^{k-1} e^{-x^k}$, for $x \geq 0$, and $k>0$, has
 $\phi(U)= - log^{1/k}(U)$. Here,   the minimal correlation, given in Equation~(\ref{mincor}), 
for different values of $k$ could only be evaluated numerically,   and is given in Table I. 
Please notice that the case  $k=1$ corresponds to $Exp(1)$ distribution that we have already discussed.  

\begin{table}[h]
\begin{center}
\begin{tabular}{|c|ccccccc|} \hline \label{weibull}
k & 4 &  3 & 2  & 1 & 0.9 & 0.8 & 0.5  \\ \hline
$\rho_\ast $ & -0.999\,&  -0.996\, &  -0.947\,&   -0.645\,&  -0.574\,&  -0.492\,& -0.193 \\ \hline
\end{tabular}\vspace{-.15in}
\end{center}
\caption{ Minimal correlation of a bivariate distribution with marginals distributed as  Weibull$(k)$, for different values of parameter $k$.}
\end{table}

\item{\bf Beta.} A random variable with $Beta(a,1)$ distribution (density $a x^{a-1}$ on $[0,1]$) can be sampled as $U^{1/a}$  and, due to symmetry, $Beta(1,b)$ can be sampled as $1-U^{1/b}$ \cite{Devroye1996}. We analyze the first case in which $\phi(U)= U^{1/a}$ and $E(\phi(U)\phi(1-U))= \mathrm{\bf B}(1/a+1,1/a+1)$, where $\mathrm{\bf B}$ stands for the beta function $\mathrm{\bf B}(x,y) = \int_0^1 t^{x-1} (1-t)^{y-1} dt$. For a special case when $a=1/n$, where $n$ is an integer, the minimum correlation can be obtained analytically by realizing  that $\mathrm{\bf B}(n+1,n+1) = \frac{(n!)^2}{(2n+1)!}.$ If we let  $m$ and $\sigma$ be the mean and standard deviation of Beta(1/n,1), then
\begin{equation*}
\rho_\ast=\frac{E(\phi(U)\phi(1-U))-m^2}{\sigma^2}= 
\end{equation*}
\\[-6ex]
\begin{equation*}
\frac{\frac{(n!)^2}{(2n+1)!}-\frac{1}{(1+n)^2}}{\frac{n^2}{(1+n)^2 (1+2n)}}=\frac{[(n+1)!]^2-(2n+1)!}{n^2 (2n)!}.
\end{equation*}

For other values of $a$,   minimal correlations can be obtained  numerically, which we  show in Table II. 
\begin{table}[h]
\begin{center}
\begin{tabular}{|c|cccccccc|} \hline  \label{tablebeta}
a & 5 &  4 & 3  & 2 & 1 & 0.8 & 0.5 & 0.3  \\ \hline 
$\rho_\ast$ & -0.795 &  -0.824 &  -0.867&  -0.931 &  -1 &  -0.989 & -0.875 &-0.634 \\ \hline
\end{tabular}\vspace{-.15in}
\end{center}
\caption{ Minimal correlation of a bivariate distribution with marginals distributed as Beta(a,1),   for different values of $a$.}
\end{table}

\item  Minimum correlations for {\bf Poisson} distribution are calculated by Shin and Pasupathy \cite{shin2009}, while {\bf log-normal} case was studied in De Veaux \cite{DeVeaux1976}, among others. We refer the readers to derivations in their papers.

\end{itemize}

\section{Multivariate algorithm} \label{multi}

In this section we propose an extension of the above algorithm to the  multivariate case. 
We start with the simplest case, where $X_1, X_2,...,X_n$ is a set of $n$ identically distributed random variables, 
each with $F$  as the marginal distribution, and identical positive pairwise correlation coefficient for each pair, $cor(X_i,X_j) = \rho^2$.  

As before, let $\phi$ be an algorithm such that $\phi(U) \sim F$, where U is a uniform random variable on [0,1]. (Although not necessary, we can set $\phi$ to equal $F^{-1}$.) Then the construction given in  Algorithm 2 below
yields a set of $n$ random variables, $(X_1,...,X_n)$, such that  $X_i \sim F$ for each $i$, and $cor(X_i,X_j)=\rho^2$ for each $i \neq j, \,\, i,j \le n$:

\begin{algorithm} [h]
\hrule
\caption{Construction of $n$ random variables, $X_1, X_2,...,X_n$, identically distributed with a 
prescribed marginal  distribution $F$ and identical positive 
pairwise correlation coefficient $\rho^2$\\[-1ex] }
\hrule \label{alg3} 
\begin{algorithmic}[1]
\STATE {\bf sample} $U,V_1,V_2 ... V_n, W_1,W_2 ... W_n \sim U(0,1)$, independently
\FOR {$i = 1 \to n$} 
\IF {$W_i < |\rho|$}
\STATE {\bf let} $X_i=\phi(U)$ \label{state12}
\ELSE
\STATE {\bf let} $X_i=\phi(V_i)$ \label{state22}
\ENDIF
\ENDFOR
\STATE RETURN $X_1,..,X_n$ \vspace{.1in}
\hrule
\end{algorithmic}
\end{algorithm}

{ Algorithm 2 will be applicable to a range of correlation values, which will depend not only on $F$
as in the bivariate case, but also on the conditions required for positive semi-definiteness or positive definiteness of
the correlation matrix. 
For example, a commonly used necessary and sufficient condition 
for positive definiteness of a matrix  is  Sylvester's criterion, which states that a matrix is positive definite if and 
only if all leading principal minors have positive determinants.  In the case of a 3-dimensional "compound symmetry" 
correlation matrix  -- a matrix where all diagonal elements are equal to 1 and all off-diagonal elements are equal to $r \in (-1,1) $ --  
Sylvester's criterion  equates to the condition that the determinant, $1 - 3r^2+ 2r^3$, is positive, or equivalently that $-0.5<r<1$.  
The matrix structure assumed by Algorithm 2 above will  thus be positive definite for any $\rho^2<1$.   The topic of conditions for 
positive semi-definiteness of a correlation matrix can be found in \cite{Devroye2010}, among others. }

The next multivariate algorithm extension is to the case when $X_1, X_2,...,X_n$ is again a set of $n$ identically distributed random variables, each with $F$  as the marginal distribution. However, now we  allow the  pairwise correlation coefficient to be different for every pair, $cor(X_i,X_j) = \rho_{ij}$, but only if each $\rho_{ij}$ can be expressed as  $\rho_i \rho_j$ where $\rho_i \in (\rho_{\ast},\rho^{\ast})$ for all $i=1,...,n$.  The construction given in  Algorithm 3 below yields a set of $n$ random variables, $(X_1,...,X_n)$, such that  $X_i \sim F$ for each $i$, and $cor(X_i,X_j)=\rho_{ij} = \rho_i \rho_j$ for each $i \neq j, \,\, i,j \le n$:

 \begin{algorithm} [h]
\hrule
\caption{Construction of $n$ random variables, $X_1, X_2,...,X_n$, identically distributed with a prescribed marginal distribution $F$, and   
pairwise correlation coefficients $\rho_{ij} = \rho_i \rho_j$\\[-1ex] } \hrule \label{alg3}
\begin{algorithmic}[1]
\STATE {\bf sample} $U,V_1,V_2 ... V_n, W_1,W_2 ... W_n  \sim U(0,1)$, independently
\FOR {$i = 1 \to n$} 
\IF {$\rho_i> 0$} 
\STATE  {\bf let} $U'=  U$ 
\ELSE
\STATE {\bf let} $U'= 1- U$ 
\ENDIF
\IF {$W_i < \rho_i / c^{\phi}(U,U')$}
\STATE {\bf let} $X_i=\phi(U')$ \label{state13}
\ELSE
\STATE {\bf let} $X_i=\phi(V_i)$ \label{state23}
\ENDIF
\ENDFOR
\STATE RETURN $X_1,..,X_n$  \vspace{.1in}
\hrule
\end{algorithmic}
\end{algorithm}

\begin{figure}[hp]
\begin{center}
{\includegraphics[width=1.8in]{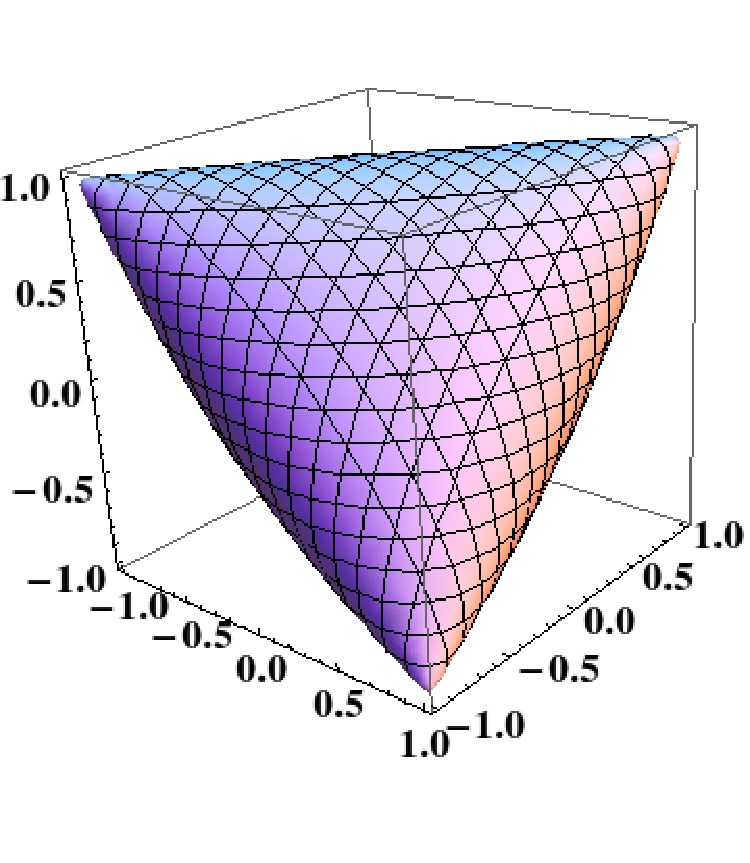}
\includegraphics[width=1.8in]{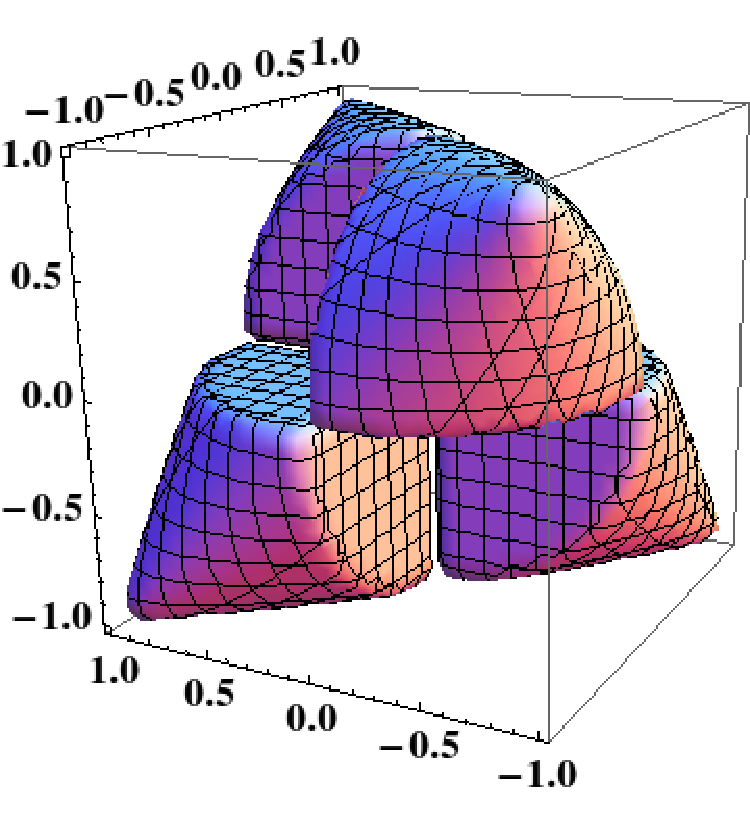}
\includegraphics[width=1.8in]{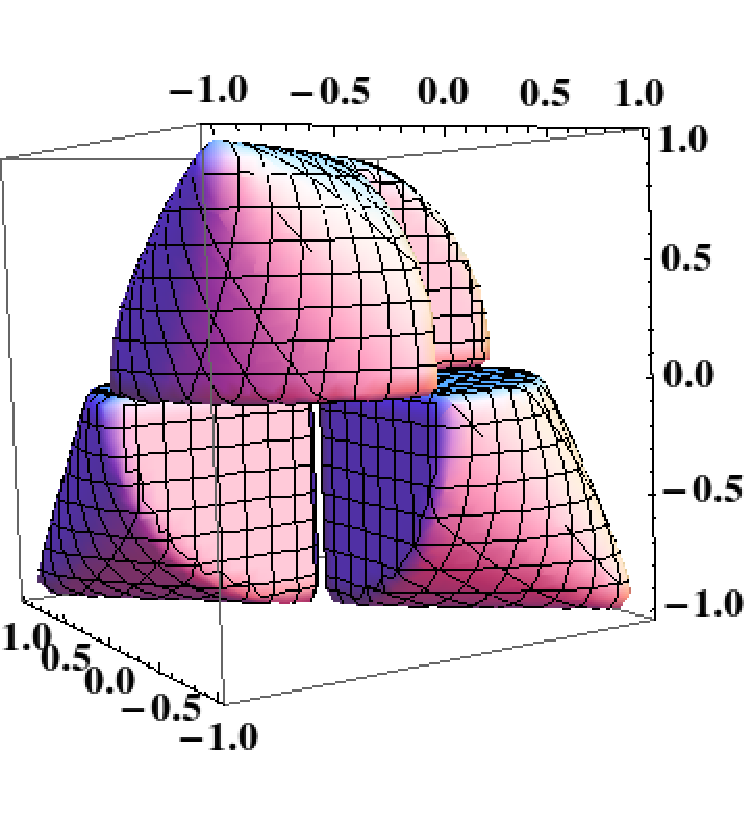}
}
\caption{(Color online) General applicable region for the 3-dimensional implementation of Algorithm 3. 
The top plot shows the full 3-dimensional domain of allowable correlation coefficients $p$, $q$, and $r$  
(shown as ranging from (-1,1) on the three coordinate axes),
which support positive definiteness of a 3-dimensional correlation matrix. The middle and bottom plots are 
alternative views of that region further restricted by the factorization requirement in Algorithm 3; 
these plots are obtained by taking the region depicted in the top plot and removing the coordinate axes and subregions where $rpq<0$.} \label{fig1}
\end{center}
\end{figure}

{ Algorithm 3  will also be applicable to a range of correlation values, which will depend not only on the choice of $F$, 
but also on the properties required of the correlation matrix. 
For example, for a general 3-dimensional correlation matrix, 
where the three correlation coefficients are $p$, $q$, and $r$, with $\left| p \right| <1,$ $\left| q \right|  <1,$ $\left| r \right|  <1,$  
Sylvester's criterion for positive definiteness equates to  the determinant being positive, $1 - p^2 - q^2 - r^2 + 2pqr>0$. }

In addition, while  Algorithm 3 allows for negative correlation between variables, the correlation coefficient factorization 
requirement imposes an added restriction.  For example,  Algorithm 3 cannot accommodate cases such as independence between $X_1$ and $X_2$, 
but dependence between $X_1$ and $X_3$ and between $X_2$ and $X_3$; nor can it accommodate correlation matrices with an odd number of negative 
correlation coefficients. This added restriction in the 3-dimensional case is shown in Figure 1. In the top plot we see  the general 
applicable 3-dimensional region for the 3 correlation coefficients $p$, $q$, and $r$ required for positive definiteness.  
Two views of the subset of that region where  Algorithm 3 is applicable  are  given the  middle and center plots. 

Only subsets of this region shown in the middle and bottom plot in Figure 1 may be applicable to specific distributions. In the case of 3-dimensional random variable with uniform marginals, the region in the middle and bottom plots in Figure 1 is fully attainable. However, in the case of a 3-dimensional 
Weibull$(0.5)$, the region shown in the middle and bottom plots in Figure 1 has to be further restricted via intersecting it with $[-0.1992,1]\times[-0.1992,1]\times[-0.1992,1]$.

Finally, we note that for any given set of correlation coefficients, $\rho_{ij}, \, i, j = 1,..., n$, the factorization into $\rho_i$ terms will generally not be unique. In particular, the factorization can be obtained by solving a set of $n \choose 2$ equations in $n$ unknowns which will, in the case of all non-zero correlations, generally yield two sets of solutions with alternate signs.  Algorithm 3 will work for any of the admissible factorizations, but choosing the factorization with the smaller number of negative $\rho_i$ coefficients is recommended, to reduce the number of comparisons  with $c^{\phi}({U},{1-U})$. In the next section we will give a concrete example of an implementation of Algorithm 3 starting with a given correlation matrix.

\subsection{Example: Algorithm for generating multivariate correlated Beta random variables}
We conclude this section with an application of Algorithm 3 for sampling a 3-dimensional random variable, $(X_1,X_2,X_3),$ with $Beta(\nu_1,\nu_2)$   marginals  and  a set of pairwise correlation coefficients $cor(X_i,X_j)=\rho_{ij}$ for  $i,j = 1,2,3.$  Among other things, a Beta density is used in practice to describe concentrations of compounds in a chemical reaction. A multivariate Beta density can thus be used to jointly describe multiple compounds, where a negative correlation would exist between a compound and its inhibitors, and a positive one between a compound and its promoters.

 Algorithm 4,  given below, is the only non-copula based algorithm for generating multivariate correlated Beta random variables for dimensions greater than 2.  This algorithm is valid for integer $\nu_1$ and $\nu_2$, and is based on generating two Gamma-distributed random variables, $G_1 \sim Gamma(\nu_1,1)$ and $G_2 \sim Gamma(\nu_2,1)$, and forming a new variable as $B = G_1 / (G_1+G_2)$, which will be distributed as $Beta(\nu_1,\nu_2)$. As $\nu_1$ and $\nu_2$ are integers, $G_1$ and $G_2$ can be obtained via a sum of $\nu_1$ and $\nu_2$ exponential random variables with mean 1, respectively.  This  example illustrates 
two facts: 1) that $\phi(\cdot)$ need not be an inverse cdf and b) that the source of randomness used in generation of a random variable  need not be a scalar. Analogous to the quantity $c^{\phi,\psi}$ defined by Equation (\ref{defc}), we will let ${\bf{U}}=(U_1,..., U_{\nu_1+\nu_2})$ and define $ \phi_{\nu_1,\nu_2} (\mathbf{U})= {\sum_1^{\nu_1} \log U_i}/{\sum_1^{\nu_1+\nu_2} \log U_i}.$

The example with 10,000 simulated 3-dimensional variables with $Beta(4,7)$  marginals,  $\rho_{12} = 0.4$, $\rho_{13} = 0.3$, and $\rho_{23} = 0.2$, resulting from Algorithm 4, is shown in Figure 2 (top). The bottom plot of Figure 2 shows an example with 10,000 simulated 3-dimensional variables with the same marginals, but with $\rho_{12} = -0.4$, $\rho_{13} = -0.3$, and $\rho_{23} = 0.3$. 
Note that for $Beta(4,7)$, $c^{\phi_{4,7}}({\bf U},{\bf 1-U}) \approx -0.71$, 
so only $\rho_{ij} \geq -0.71$ can be considered for generating $Beta(4,7)$ using Algorithm 4.

\begin{figure}[htp]
\begin{center}
{\includegraphics[height=2.5in,angle=0]{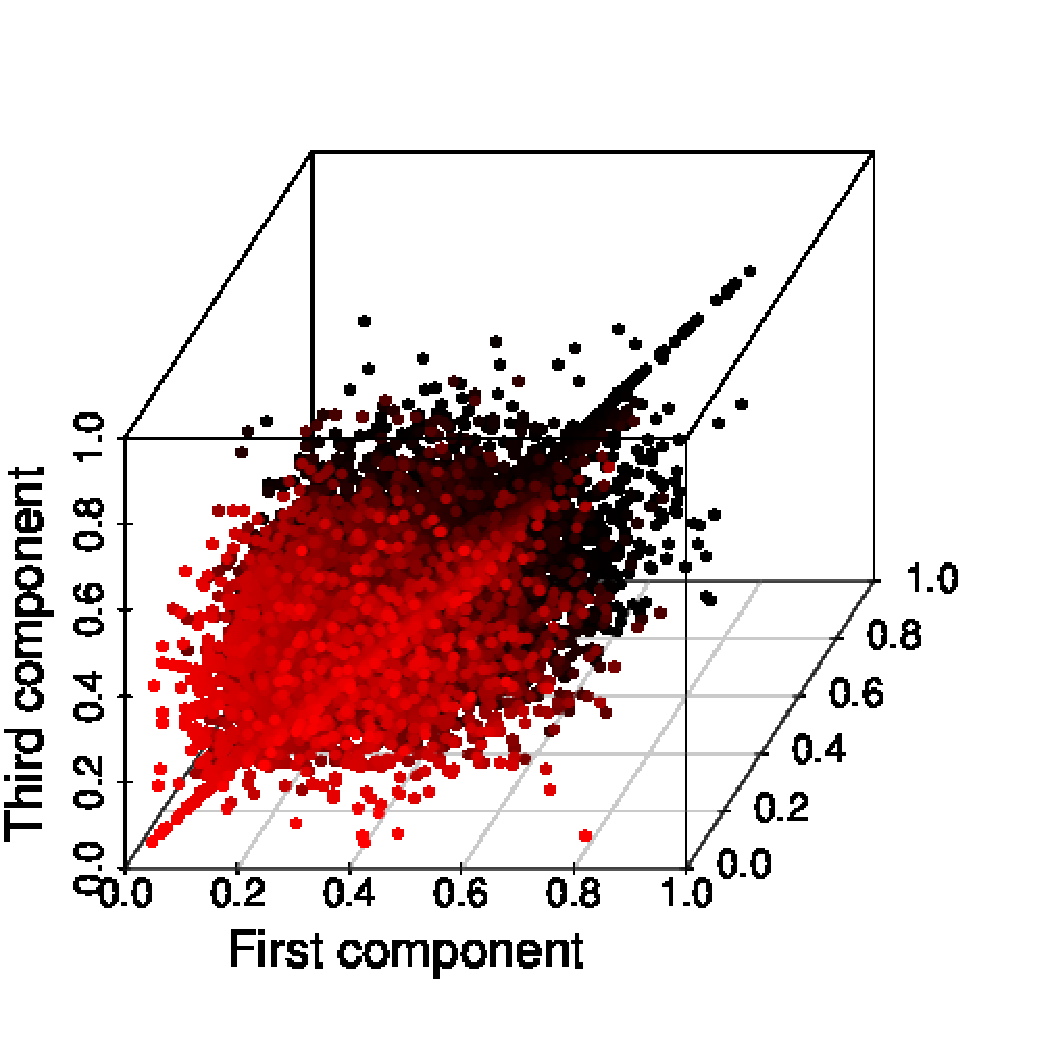}
\includegraphics[height=2.5in,angle=0]{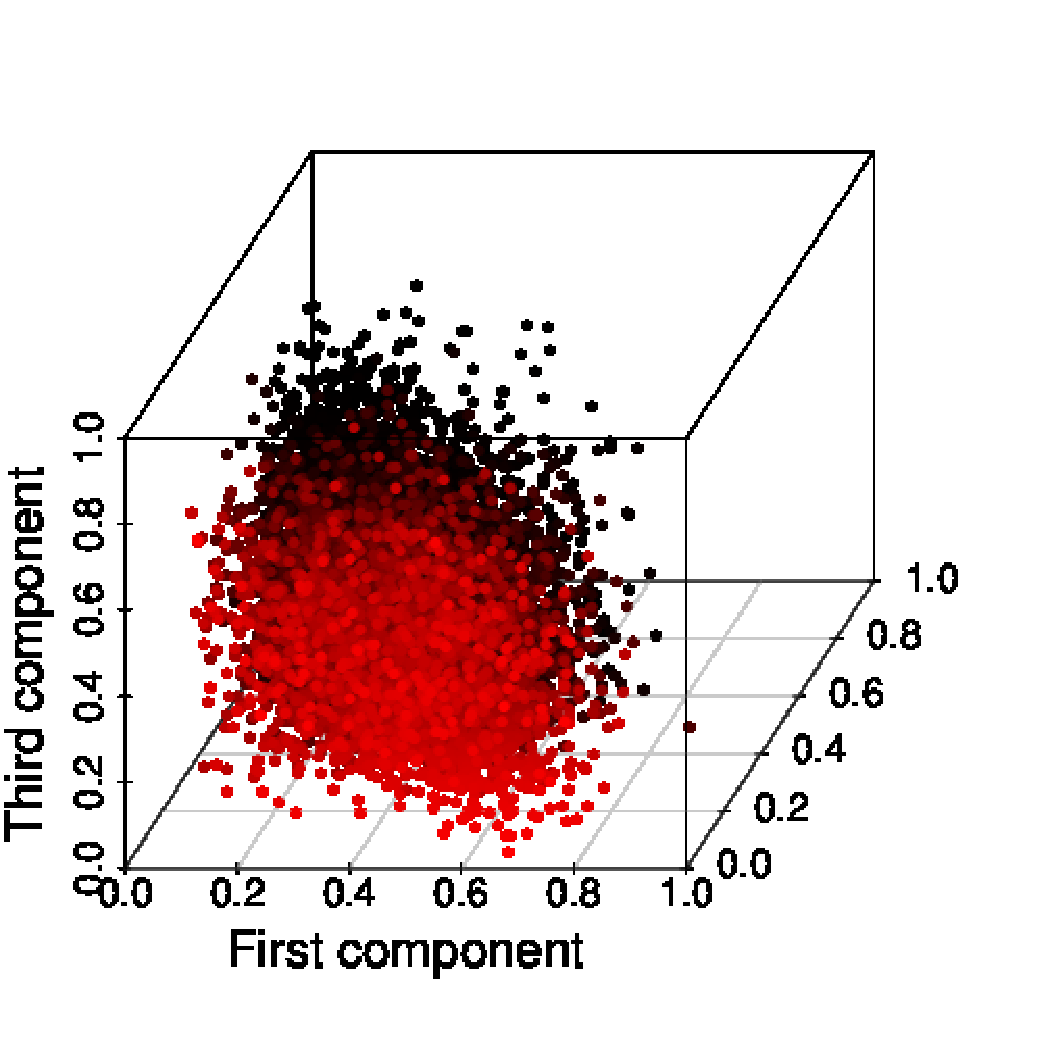}
}
\caption{(Color online) An example with 10,000 simulated 3-dimensional $Beta(4,7)$ variables resulting from  Algorithm 4. Top panel: $\rho_{12} = 0.4$, $\rho_{13} = 0.3$, $\rho_{23} = 0.2$; Bottom panel: $\rho_{12} = -0.4$, $\rho_{13} = -0.3$, and $\rho_{23} = 0.3$.} \label{fig2}
\end{center}
\end{figure}

\section{Conclusions}\label{Conclusions}
The algorithm presented in this paper is a simple generalization of the trivariate-reduction method for generation of  
multivariate samples with specified  marginal distributions and correlation matrix. 
{ In comparison with the copulas it is  
 simpler in that it is based only on marginal distributions and a correlation matrix 
and does not require a whole multivariate distribution specification. On the other hand
it  is  exact and more transparent to implement than copulas. Additionally, we generate samples directly from uniform 
random variables, the immediate output from random number generators, which may 
be more desirable and faster than going through others distributions, such as Gaussians, as in 
many other methods. } 

The algorithm is applicable to all distributions with finite variances, and, in 
the bivariate case, can accommodate the entire range of theoretically feasible correlations. Its 
major computational difficulty is related to determination of exact pairwise correlation ranges, a question of 
theoretical and practical value  {\it per se}, which has to be answered {\em once} for every set 
of marginal distributions.  { We emphasize  that lower and upper bounds for the correlation coefficient
  actually depend on the family of marginal distributions in question, and that the commonly used $[-1,1]$
interval can be inappropriate in many applications.  In this paper we have 
presented exact ranges for some common distributional examples so that the implementation of the algorithms is straightforward. }

\begin{algorithm} [ht]
\hrule
\caption{Construction of 3-dimensional  $Beta(\nu_1,\nu_2)$ random variable, with $Beta(\nu_1,\nu_2)$  marginals and with  correlation coefficients
 $\rho_{12}$, $\rho_{13}$, $\rho_{23}$.\\[-1ex]
} \hrule \label{alg5}
\begin{algorithmic}[1]
\STATE {\bf sample} ${\bf U} = \left\{ U_1^1,..,U_1^{\nu_1}, U_2^1,..,U_2^{\nu_2}\right\},  W_1,W_2, W_3  \sim U(0,1)$, independently
\IF { { any of the eigenvalues of the given correlation matrix are negative} }
\STATE {\bf stop: matrix not positive {semi-definite}}
\ELSE
\IF { $\rho_{12}\rho_{13}\rho_{23} < 0$, only one of $\rho_{12}$, $\rho_{13}$, $\rho_{23}$ is 0, or 
any $\rho_{ij} \leq c^{\phi_{\nu_1,\nu_2}}({\bf U},{\bf 1-U})$}
\STATE {\bf stop: algorithm not applicable.}
\ELSE
\IF {  $\rho_{12} = \rho_{13} = \rho_{23} = 0$}
\STATE {\bf let} $\rho_1=\rho_2=\rho_3=0$
\ELSE
\IF {  $\rho_{ij}=0$, $\rho_{ik}=0$ and $\rho_{jk}\neq 0$}  
\STATE {\bf let} $\rho_i = 0$, $\rho_j = 1$ and $\rho_k = \rho_{ik}$
\ELSE
\STATE {\bf let} $\rho_2 = \sqrt{\rho_{12} \rho_{23}/\rho_{13}}$, $\rho_1 = \rho_{12}/\rho_2$, $\rho_3 = \rho_{23}/\rho_2$
\IF {  $\rho_i \leq c^{\phi_{\nu_1,\nu_2}}({\bf U},{\bf 1-U})$ for any $i$}
\STATE {\bf warning: algorithm will produce only approximate results for negative correlations.}
\ELSE
\FOR {$i = 1 \to 3$} 
\IF {$\rho_i> 0$} 
\STATE  {\bf let} ${\bf U'}=  {\bf U}$ 
\ELSE
\STATE {\bf let} ${\bf U'}= {\bf 1}- {\bf U}$
\ENDIF
\IF {$W_i < \rho_i/ c^{\phi_{\nu_1,\nu_2}}({\bf U},{\bf U'})$}
\STATE {\bf let} $G_1 = {\sum_{j=1}^{\nu_1} -\log(U_1^{'j})}$
\STATE {\bf let} $G_2 = {\sum_{j=1}^{\nu_2} -\log(U_2^{'j})}$
\ELSE
\STATE {\bf sample} $V_{1}^1,..,V_1^{\nu_1}, V_{2}^1,..,V_{2}^{\nu_2}  \sim U(0,1)$, independently
\STATE {\bf let} $G_1 = {\sum_{j=1}^{\nu_1} -\log(V_{1}^j)}$
\STATE {\bf let} $G_2 = {\sum_{j=1}^{\nu_2} -\log(V_{2}^j)}$
\ENDIF
\STATE {\bf let} $X_i=G_1/(G_1+G_2) $ \label{state2}
\ENDFOR
\ENDIF
\ENDIF
\ENDIF
\ENDIF
\ENDIF
\STATE RETURN $X_1,X_2,X_3$ 
\vspace{.1in}
\hrule
\end{algorithmic}
\end{algorithm}
\clearpage

\noindent {\bf Acknowledgements}
We thank Fabio Machado for introducing this problem to us and to Mark Huber, Xiao-Li Meng, and David Bortz for insightful discussions. This work was supported in part by  NIGMS (U01GM087729), NIH (R21DA027624), and NSF (DMS-1007823).

\section*{Appendix}

To prove the Exponential distribution result from Section~\ref{examples}:
\begin{eqnarray*}
\int_0^1 \log(x) \log(1-x) dx = 2 -\frac{1}{6} \pi^2,
\end{eqnarray*}
\noindent 
we use a representation of $\log(x)$ as a Maclaurin series for $x \in (0,1)$:
\begin{eqnarray*}
&\log(x)\log(1-x)=-\sum_{i=1}^{\infty} \frac{x^i}{i} \sum_{j=1}^{\infty} (-1)^{j+1}\frac{(x-1)^j}{j}=& \\
& \sum_{i=1}^{\infty}  \sum_{j=1}^{\infty} \frac{x^i}{i}\frac{(1-x)^j}{j}.&
\end{eqnarray*}
Observe that $\lim_{x\to0} \log(x)\log(1-x)=\lim_{x\to1} \log(x)\log(1-x)=0$, and the double sum equals $\log(x)\log(1-x)$ for all $x \in [0,1]$.
Furthermore, 
\begin{eqnarray} 
&&\int_0^1 \log(x)\log(1-x)dx=\sum_{i=1}^{\infty}  \sum_{j=1}^{\infty} \frac{1}{ij} \beta (i+1,j+1)
 \nonumber \\
&&= \sum_{i=1}^{\infty}  \sum_{j=1}^{\infty} \frac{1}{ij} \frac{i! j!}{(i+j+1)!} \nonumber \\
&&=\sum_{i=1}^{\infty}  \sum_{j=1}^{\infty}  \frac{(i-1)! (j-1)!}{(i+j+1)!} =\sum_{i=1}^{\infty}  \sum_{j=1}^{\infty}  \frac{(i-1)!}{(i+2)!{i+j+1 \choose j-1}}\nonumber \\
&&=\sum_{i=1}^{\infty} \frac{1}{i(i+1)(i+2)}  \sum_{j=1}^{\infty}  \frac{1}{{i+j+1 \choose j-1},} \\[4ex]
&&  \nonumber \label{zadnja}
\end{eqnarray}
where $\beta(i+1,j+1) = \int_0^1x^i(1-x)^{j}dx= \frac{i!j!}{(i+j+1)!}$ is a standard presentation of beta function with $i,j$ integers. 

To proceed from here we use the Corollary 3.7 in \cite{sury}:
\begin{eqnarray*}
\sum_{k=0}^{\infty}  \frac{1}{{n+k \choose k}}= \frac{n}{n-1}
\end{eqnarray*}
so the last $j$-sum in (\ref{zadnja}) equals $(i+2)/(i+1)$ and
\begin{eqnarray}\label{sum}
\int_0^1 \log(x)\log(1-x)dx&=& \sum_{i=1}^{\infty} \frac{1}{i(i+1)^2} .
\end{eqnarray}
To prove that the above series converges to  $2-\pi^2/6$ recall that $\sum_{i=1}^{\infty} 1/i^2 =\pi^2/6$. Now we add that series to (\ref{sum}) and show that it adds up to 2:
\begin{eqnarray*}
 \sum_{i=1}^{\infty} \frac{1}{i(i+1)^2} + \frac{\pi^2}{6}= \sum_{i=1}^{\infty} \frac{1}{i(i+1)^2} + \sum_{i=1}^{\infty} \frac{1}{i^2}\\
 =1+\sum_{i=1}^{\infty} \frac{1}{(i+1)^2}(\frac{1}{i} +1)  \\
 =1+\sum_{i=1}^{\infty} \frac{1}{i(i+1)}= 1+\sum_{i=1}^{\infty}( \frac{1}{i}-\frac{1}{i+1})=1+1=2.  \end{eqnarray*}

\bibliography{FrechetAll}
\bibliographystyle{apsrev}

\end{document}